\documentclass{amsart}[12pt]

\theoremstyle{definition}

\usepackage{amscd,amssymb}

\begin{document}

\title[]
{Weak polynomial identities for a vector space\\
with a symmetric bilinear form}
\author[]{Vesselin S. Drensky, Plamen E. Koshlukov}

\address{Permanent addresses in 2019:}
\address{Vesselin Drensky: Institute of Mathematics and Informatics,
Bulgarian Academy of Sciences,
Acad. G. Bonchev Str., Block 8, 1113 Sofia, Bulgaria}
\email{drensky@math.bas.bg}

\address{Plamen Koshlukov: Department of Mathematics, IMECC, UNICAMP,
S\'ergio Buarque de Holanda 651, Campinas,
SP 13083-859, Brazil}
\email{plamen@ime.unicamp.br}

\thanks{Published as Vesselin S. Drensky, Plamen E. Koshlukov, Weak polynomial identities for a vector space with a symmetric bilinear form (Bulgarian summary),
Mathematics and Mathematical Education, 1987. Proceedings of the Sixteenth Spring Conference
of the Union of Bulgarian Mathematicians, Sunny Beach, April 6-10, 1987,
213-219, Publishing House of the Bulgarian Academy of Sciences, Sofia, 1987. Zbl 0658.16012, MR949941 (89j:15043).}
\date{}

\maketitle

\begin{abstract}
Let $V_k$ be a $k$-dimensional vector space with a non-degenerate symmetric bilinear form and let $C_k$ be the Clifford algebra on $V_k$.
The weak polynomial identities of the pair $(C_k,V_k)$ are investigated. It is proved that all they follow from $[x_1^2,x_2]=0$
when $k=\infty$ and from $[x_1^2,x_2]=0$ and $S_{k+1}(x_1,\ldots,x_{k+1})=0$ when $k<\infty$. The Specht property of the weak identity $[x_1^2,x_2]=0$
is established as well.
\end{abstract}

\section{Introduction}
The weak polynomial identities have been introduced by Razmyslov \cite{R} in his studying of the $2\times 2$ matrix algebra. They can be successfully
applied in the investigation of the polynomial identities for algebras close to associative \cite{V}, \cite{I}.

All our considerations are made over a fixed field $K$ of characteristic 0.
Let $R$ be an associative algebra and let $G$ be a vector subspace of $R$ generating $R$
as an algebra. The polynomial $f(x_1,\ldots,x_n)$ of the free associative algebra $K\langle X\rangle$ is called a weak identity for the pair $(R,G)$ if
$f(g_1,\ldots,g_n)=0$ for all $g_1,\ldots,g_n\in G$. The weak identities of $(R,G)$ form an ideal $T(R,G)$ of $K\langle X\rangle$.

There are various possibilities for defining consequences of a weak identity depending on the properties of $G$.
Generally let $\Omega\subset K\langle X\rangle$ be a family
of polynomials such that $\omega(g_1,\ldots,g_n)\in G$ for every $\omega(x_1,\ldots,x_n)\in\Omega$ and for all $g_1,\ldots,g_n\in G$.
The weak identity $f_2(x_1,\ldots,x_m)$
is an $\Omega$-consequence of $f_1(x_1,\ldots,x_n)$ if $f_2$ belongs to the ideal of $K\langle X\rangle$
generated by $\{f_1(\omega_1,\ldots,w_n)\mid \omega_i\in\Omega,i=1,\ldots,n\}$.

{\bf Examples:} 1. When $G=R$ and $\Omega=K\langle X\rangle$ we obtain the ordinary polynomial identities
with the usual rules for finding consequences of an identity.

2. Let $[G,G]\subseteq G$, i.e., let $G\subset R$ be a subalgebra of the adjoint Lie algebra $R^{(-)}$ of $R$
and let $\Omega=L(X)\subset K\langle X\rangle$ be the free Lie algebra
canonically embedded into $K\langle X\rangle$. Then the weak Lie identities are obtained \cite{R}.

3. Assume $G$ is closed with respect to the operation $\displaystyle g\circ h=\frac{1}{2}(gh+hg)$,
i.e., $G$ is a subalgebra of the Jordan algebra $R^{(+)}$. In this case
we obtain the weak Jordan identities.

4. Let $\Omega=\text{sp}(X)$ be the vector subspace of $K\langle X\rangle$ spanned by $X$.
The ideal $T(R,G)$ of the weak identities for every pair $(R,G)$ is $\Omega$-closed.
Hence if $f(x_1,\ldots,x_n)\in T(R,G)$ then $\displaystyle f\left(\sum\alpha_{1j}x_j,\ldots,\sum\alpha_{nj}x_j\right)\in T(R,G)$ for every $\alpha_{ij}\in K$.
The free algebra $K\langle X\rangle$ is isomorphic to the tensor algebra of $\text{sp}(X)$
and the general linear group $GL=GL(\text{sp}(X))$ acts canonically on $K\langle X\rangle$.
Obviously $T(R,G)$ is a $GL$-invariant ideal. Assume $H=\{h_i(x_1,\ldots,x_{n_i})\mid i\in I\}$ is a subset of $K\langle X\rangle$.
The ideal generated by $\{h_i(\omega_1,\ldots,\omega_{n_i})\mid i\in I,\omega_j\in\text{sp}(X)\}$
is called a $GL$-ideal generated by $H$ and its elements are $GL$-consequences of $H$.

Let $G=V_k$ be a $k$-dimensional vector space with a non-degenerate symmetric bilinear form $\langle\text{ },\text{ }\rangle$,
and let $R=C_k$ be the Clifford algebra of $V_k$.
Furthermore we set $V=V_{\infty}$ and $C=C_{\infty}$.

The investigation of the weak identities for the pair $(C_k,V_k)$ has been initiated by Il'tyakov \cite{I}.
Our main purpose is to find a basis for the identities of $(C,V)$.

{\bf Theorem 1.} {\it The $GL$-ideal $T(C,V)$ of all weak identities in the pair $(C,V)$ is generated by $[x_1^2,x_2]$.}

The proof of this result uses the technique of the representations of the general linear group and it is in the spirit of \cite{D}.

{\bf Corollary 1.} {\it The weak identities for the pair $(C_k,V_k)$ are $GL$-consequences of
\begin{equation}
[x_1^2,x_2]
\end{equation}
and of the standard polynomial}
\[
S_{k+1}(x_1,\ldots,x_{k+1})=\sum(-1)^{\sigma}x_{\sigma(1)}\cdots x_{\sigma(k+1)}.
\]

Assume $M_2$ is the $2\times 2$ matrix algebra and $sl_2$ is its Lie subalgebra of all traceless matrices.
Razmyslov \cite{R} found a basis for the weak Lie identities of $(M_2,sl_2)$.
As a consequence of Theorem 1 we obtain

{\bf Corollary 2.} (i) {\it The $GL$-ideal $T(M_2,sl_2)$ is generated by
\[
[x_1^2,x_2]\text{ and }S_4(x_1,x_2,x_3,x_4).
\]
}

(ii) \cite{R} {\it The weak Lie identities of $(M_2,sl_2)$ follow from} (1).

Applying ideas from \cite{K} we prove that the weak identity (1) satisfies the Specht property.

{\bf Theorem 2.} {\it Every $GL$-ideal of $K\langle X\rangle$ containing the polynomial $[x_1^2,x_2]$ is finitely generated.}

The results of this note are very close related to the problems for finding the bases of the ordinary and the weak identities for the Jordan algebra
of a symmetric bilinear form. We hope they represent a step to the solution of these important problems in the theory of Jordan algebras
with polynomial identities.

\section{Identities in the pair $(C,V)$}

Let us denote by $U$ the $GL$-ideal of $K\langle X\rangle$ generated by the polynomial $[x_1^2,x_2]$, $F=K\langle X\rangle/U$, and
$F_m=K\langle x_1,\ldots,x_m\rangle/(K\langle x_1,\ldots,x_m\rangle\cap U)$. In the sequel we shall work in the algebras $F$ and $F_m$.
We shall also use other letters, e.g., $y_1,y_2,\ldots$ for the generators of $F$.

{\bf Lemma 1.} {\it For each $n>1$ an equality of the form
\begin{equation}
x_1y_1\cdots y_nx_2-x_2y_1\cdots y_nx_1=\sum_iA_{n_i}[x_1,x_2]B_{n_i}
\end{equation}
holds for $F$. Here $A_{n_i}$ and $B_{n_i}$ are homogeneous polynomials in the variables $y_1,\ldots,y_n$ and $\deg A_{n_i}>0$.
In particular for $n=2$
\begin{equation}
x_1y_1y_2x_2-x_2y_1y_2x_1=\frac{1}{2}\{(y_1y_2+y_2y_1)[x_1,x_2]+y_1[x_1,x_2]y_2-y_2[x_1,x_2]y_1\}
\end{equation}
holds.}

\begin{proof}
The linearization of (1) is
\[
f_1(x_1,x_2,y)=[x_1\circ x_2,y]=0.
\]
Hence we obtain
\[
0=f_1(x_1,x_2,y)+2f_1(x_1,y,x_2)=x_1yx_2-x_2yx_1+[x_1,x_2]\circ y,
\]
i.e.,
\begin{equation}
x_1yx_2-x_2yx_1=-[x_1,x_2]\circ y.
\end{equation}
Bearing in mind that
\[
f_2(x_1,x_2,y_1,\ldots,y_n)=[x_1\circ y_1,y_2]y_3\cdots y_nx_2=0
\]
we deduce
\[
2\{f_2(x_1,x_2,y_1,\ldots,y_n)-f_2(x_2,x_1,y_1,\ldots,y_n)\}
=x_1y_1\cdots y_nx_2-x_2y_1\cdots y_nx_1
\]
\[
+y_1(x_1y_2\cdots y_nx_2-x_2y_2\cdots y_nx_1)-y_2(x_1y_1y_3\cdots y_nx_2-x_2y_1y_3\cdots y_nx_1)
\]
\[
-y_1y_2(x_1y_3\cdots y_nx_2-x_2y_3\cdots y_nx_1)=0.
\]
Hence we can express $x_1y_1\cdots y_nx_2-x_2y_1\cdots y_nx_1$ as a linear combination of polynomials beginning with $y_1$ or $y_2$
and skew-symmetric in $x_1$ and $x_2$. Using (4) we establish (3). The equality (2) follows by obvious induction on $n$.
\end{proof}

{\bf Lemma 2.} {\it Weak identities of the following form are $GL$-consequences of} (1):
\begin{equation}
\begin{split}
\sum(-1)^{\sigma}x_{\sigma(1)}\cdots x_{\sigma(k)}yx_{\sigma(k+1)}\cdots x_{\sigma(n)}\\
=\alpha_{nk}yS_n(x_1,\ldots,x_n)+\beta_{nk}S_n(x_1,\ldots,x_n)y,
\end{split}
\end{equation}
{\it where $\alpha_{nk},\beta_{nk}\in K$, $\alpha_{nk}=\beta_{n,n-k}$, $k=1,\ldots,n-1$, furthermore}
\[
\alpha_{n1}=-\frac{n-1}{n},\beta_{n1}=\frac{(-1)^{n-1}}{n};
\]
\begin{equation}
x_iS_n(x_1,\ldots,x_n)=(-1)^{n-1}S_n(x_1,\ldots,x_n)x_i\text{ for }i=1,\ldots,n.
\end{equation}

\begin{proof}
(i) An involution $\ast$ is defined in the algebra $K\langle X\rangle$ by
\[
\left(\sum\alpha_ix_{i_1}\cdots x_{i_n}\right)^{\ast}=\sum\alpha_ix_{i_n}\cdots x_{i_1}.
\]
Since $[x_1^2,x_2]^{\ast}=-[x_1^2,x_2]$, the action of $\ast$ is carried over to $F$. Thus if we apply $\ast$ to (5) we obtain
\[
\sum(-1)^{\sigma}x_{\sigma(n)}\cdots x_{\sigma(k+1)}yx_{\sigma(k)}\cdots x_{\sigma(1)}
\]
\[
=\sum(-1)^{\sigma}\{\alpha_{nk}x_{\sigma(n)}\cdots x_{\sigma(1)}y+\beta_{nk}yx_{\sigma(n)}\cdots x_{\sigma(1)}\}.
\]
Replacing the positions of $\sigma(n)$ and $\sigma(1)$, $\sigma(n-1)$ and $\sigma(2)$, etc. we change simultaneously the sign of all summands.
Hence we obtain that $\alpha_{nk}=\beta_{n,n-k}$ holds in (5).

First we shall prove (5) by induction on $n$ for $k=1$ (and therefore for $k=n-1$) and then we shall examine the case of an arbitrary $k$.
For $n=2$ and $k=1$ the equality (4) gives
\[
\sum(-1)^{\sigma}x_{\sigma(1)}yx_{\sigma(2)}=-\frac{1}{2}yS_2(x_1,x_2)-\frac{1}{2}S_2(x_1,x_2)y.
\]
For $n=3$, $k=1$ we obtain from (3)
\[
\sum(-1)^{\sigma}x_{\sigma(1)}(yx_{\sigma(2)})x_{\sigma(3)}
\]
\[
=\frac{1}{2}\sum(-1)^{\sigma}\{(yx_{\sigma(2)}+x_{\sigma(2)}y)x_{\sigma(1)}x_{\sigma(3)}
+yx_{\sigma(1)}x_{\sigma(2)}x_{\sigma(3)}-x_{\sigma(1)}x_{\sigma(2)}x_{\sigma(3)}y\}
\]
\[
=-\frac{1}{2}\sum(-1)^{\sigma}x_{\sigma(1)}yx_{\sigma(2)}x_{\sigma(3)}-yS_3(x_1,x_2,x_3)+\frac{1}{2}S_3(x_1,x_2,x_3)y.
\]
Therefore it holds
\[
\sum(-1)^{\sigma}x_{\sigma(1)}yx_{\sigma(2)}x_{\sigma(3)}=-\frac{2}{3}yS_3(x_1,x_2,x_3)+\frac{1}{3}S_3(x_1,x_2,x_3)y.
\]
Let $n>3$, $k=1$. Then
\[
\sum(-1)^{\sigma}(x_{\sigma(1)}yx_{\sigma(2)}\cdots x_{\sigma(n-1)})x_{\sigma(n)}
\]
\[
=\sum(-1)^{\sigma}\left\{-\frac{n-2}{n-1}yx_{\sigma(1)}\cdots x_{\sigma(n)}+\frac{(-1)^{n-2}}{n-1}x_{\sigma(1)}
(x_{\sigma(2)}\cdots x_{\sigma(n-1)}yx_{\sigma(n)})\right\}
\]
\[
=-\frac{n-2}{n-1}yS_n(x_1,\ldots,x_n)
\]
\[
+\frac{(-1)^{n-2}}{n-1}\sum(-1)^{\sigma}\left\{\frac{(-1)^{n-2}}{n-1}x_{\sigma(1)}yx_{\sigma(2)}\cdots x_{\sigma(n)}
-\frac{n-2}{n-1}x_{\sigma(1)}\cdots x_{\sigma(n)}y\right\}
\]
and
\[
\left(1-\frac{1}{(n-1)^2}\right)\sum(-1)^{\sigma}x_{\sigma(1)}yx_{\sigma(2)}\cdots x_{\sigma(n)}
\]
\[
=-\frac{n-2}{n-1}yS_n(x_1,\ldots,x_n)-\frac{(-1)^{n-2}(n-2)}{(n-1)^2}S_n(x_1,\ldots,x_n)y,
\]
\[
\sum(-1)^{\sigma}x_{\sigma(1)}yx_{\sigma(2)}\cdots x_{\sigma(n)}=-\frac{n-1}{n}yS_n(x_1,\ldots,x_n)+\frac{(-1)^{n-1}}{n}S_n(x_1,\ldots,x_n)y.
\]
Now we assume $1<k<n$. We shall suppose that (5) holds for all smaller values of $n$. Then we have
\[
\sum(-1)^{\sigma}x_{\sigma(1)}(x_{\sigma(2)}\cdots x_{\sigma(k)}yx_{\sigma(k+1)}\cdots x_{\sigma(n)})
\]
\[
=\sum(-1)^{\sigma}\alpha_{n-1,k-1}x_{\sigma(1)}yx_{\sigma(2)}\cdots x_{\sigma(n)}
+\beta_{n-1,k-1}S_n(x_1,\ldots,x_n)y
\]
\[
=(\alpha_{n-1,k-1}\alpha_{n1}yS_n(x_1,\ldots,x_n)
+(\alpha_{n-1,k-1}\beta_{n1}+\beta_{n-1,k-1})S_n(x_1,\ldots,x_n)y.
\]

(ii) It suffices to prove the equality (6) for $i=1$ only. Obviously
\[
x_1S_n(x_1,\ldots,x_n)=x_1^2\sum_{\sigma(1)=1}(-1)^{\sigma}x_{\sigma(2)}\cdots x_{\sigma(n)}
+x_1\sum_{\sigma(1)\not=1}(-1)^{\sigma}x_{\sigma(1)}\cdots x_{\sigma(n)}.
\]
Similarly,
\[
S_n(x_1,\ldots,x_n)x_1=\sum_{\tau(n)=1}(-1)^{\tau}x_{\tau(1)}\cdots x_{\tau(n-1)}x_1^2
+\sum_{\tau(n)\not=1}(-1)^{\tau}x_{\tau(1)}\cdots x_{\tau(n)}x_1.
\]
By (1) $x_1^2$ lies in the center of $F$, hence
\[
(-1)^{\sigma}x_1^2x_{\sigma(2)}\cdots x_{\sigma(n)}=(-1)^{\sigma}x_{\sigma(2)}\cdots x_{\sigma(n)}x_1^2
=(-1)^{n-1}(-1)^{\tau}x_{\tau(1)}\cdots x_{\tau(n-1)}x_1^2,
\]
where $\sigma(1)=\tau(n)=1$, $\tau(1)=\sigma(2),\ldots,\tau(n-1)=\sigma(n)$. Now it follows from (5):
\[
x_1\sum(-1)^{\sigma}x_{\sigma(1)}\cdots x_{\sigma(k-1)}x_1x_{\sigma(k+1)}\cdots x_{\sigma(n)}
\]
\[
=(-1)^kx_1\{\alpha_{n-1,k}x_1S_{n-1}(x_2,\ldots,x_n)+\beta_{n-1,k}S_{n-1}(x_2,\ldots,x_n)x_1\},
\]
\[
\sum(-1)^{\tau}x_{\tau(1)}\cdots x_{\tau(n-k-1)}x_1x_{\tau(n-k+1)}\cdots x_{\tau(n)}x_1
\]
\[
=(-1)^{n-k-1}\{\alpha_{n-1,n-1-k}x_1S_{n-1}(x_2,\ldots,x_n)x_1+\beta_{n-1,n-1-k}S_{n-1}(x_2,\ldots,x_n)x_1^2\}
\]
and we establish (6) from the conditions $\alpha_{n-1,n-1-k}=\beta_{n-1,k}$ and $\beta_{n-1,n-1-k}=\alpha_{n-1,k}$.
\end{proof}

{\bf Lemma 3.} {\rm (i)} {\it Let $Y_1,\ldots,Y_{n-1}$ be monomials in $y_1,y_2,\ldots$. Then there exist monomials $D_{n_i},E_{n_i}$ in
$y_1,y_2,\ldots$ such that it holds in $F$}
\begin{equation}
\sum(-1)^{\sigma}x_{\sigma(1)}Y_1x_{\sigma(2)}Y_2\cdots Y_{n-1}x_{\sigma(n)}=\sum_i D_{n_i}S_n(x_1,\ldots,x_n)E_{n_i}.
\end{equation}

{\rm (ii)} {\it Let $Y_1,\ldots,Y_{n-1}$ be monomials in $x_1,\ldots,x_n$. Then there exists a polynomial $D(x_1,\ldots,x_n)$ such that in $F$}
\begin{equation}
\sum(-1)^{\sigma}x_{\sigma(1)}Y_1x_{\sigma(2)}Y_2\cdots Y_{n-1}x_{\sigma(n)}=S_n(x_1,\ldots,x_n)\cdot D(x_1,\ldots,x_n).
\end{equation}

\begin{proof}
(i) We shall carry out an induction on $n$ and on $d=\deg Y_1+\cdots+\deg Y_{n-1}$.
The basis for the induction $n=2$ follows from (2). For each $n$ and for $d=1$ (7) is a consequence of (5).
We assume that (7) holds for $n-1$ and for all values of $d$. When $d>1$ we obtain
\[
\begin{split}
\sum(-1)^{\sigma}x_{\sigma(1)}Y_1(x_{\sigma(2)}Y_2\cdots Y_{n-1}x_{\sigma(n)})\\
=\sum_{i,\sigma}(-1)^{\sigma}x_{\sigma(1)}Y_1D_{n-1,i}x_{\sigma(2)}\cdots x_{\sigma(n)}E_{n-1,i}.
\end{split}
\]
We apply the induction to the summands with $\deg E_{n-1,i}>0$. When $\deg E_{n-1,i}=0$ we use (2) for
$x_{\sigma(1)}Y_1D_{n-1,i}x_{\sigma(2)}-x_{\sigma(2)}Y_1D_{n-1,i}x_{\sigma(1)}$.
Then some variable in $Y_1D_{n-1,i}$ appears on the left-hand side of the sum and we can apply again the induction on $d$.

(ii) The equality (8) is a direct consequence of (6) and (7).
\end{proof}

{\it Proof of Theorem 1.}
First we note that the pair $(C,V)$ satisfies the identity (1) since for each vector $v\in V$, $v^2=\langle v,v\rangle\in K$ and $v^2$ lies in
the center of $C$. Now let $\lambda=(\lambda_1,\ldots,\lambda_r)$, $\lambda_1\geq\cdots\geq\lambda_r>0$, be a partition and let the columns
of the diagram $[\lambda]$ have lengths $r=r_1,r_2,\ldots,r_p$, respectively. We choose an integer $m\geq r$. Since the ideal $U$ of $K\langle X\rangle$
is $GL$-invariant, the algebra $F_m$ is a $GL_m$-module. Assume $N_m(\lambda)$ is an irreducible $GL_m$-submodule of $F_m$ corresponding to
the partition $\lambda$. It is known that $N_m(\lambda)$ is generated by a multihomogeneous polynomial $f_{\lambda}(x_1,\ldots,x_r)$ of
degree $\lambda_i$ in $x_i$. Furthermore
\[
f_{\lambda}(x_1,\ldots,x_r)=\sum_Y\sum_{\sigma}(-1)^{\sigma}Y_0x_{\sigma(1)}Y_1x_{\sigma(2)}\cdots Y_{r-1}x_{\sigma(r)}Y_r,
\]
where $Y_0,Y_1,\ldots,Y_r$ are monomials in $x_1,\ldots,x_r$. By (8)
\[
f_{\lambda}(x_1,\ldots,x_r)=S_r(x_1,\ldots,x_r)\cdot D(x_1,\ldots,x_r)
\]
holds in $F_m$ and the polynomial $D(x_1,\ldots,x_r)$ generates the irreducible $GL_m$-module $N_m(\lambda_1-1,\lambda_2-1,\ldots,\lambda_r-1)$.

Applying several times the equality (8) we obtain that in $F_m$
\[
f_{\lambda}(x_1,\ldots,x_r)=\alpha\prod S_{r_i}(x_1,\ldots,x_{r_i}),\quad \alpha\in K.
\]
Therefore all isomorphic irreducible submodules of $F_m$ are glued together. Hence $F_m$ is a submodule of $\sum N_m(\mu)$ where the summation is
over all partitions $\mu=(\mu_1,\ldots,\mu_r)$, $\mu_1\geq\cdots\geq\mu_r\geq 0$. Assume $e_1,e_2,\ldots$ is a basis of $V$ such that
\[
\langle e_i,e_i\rangle\not=0,\langle e_i,e_j\rangle=0, i\not= j.
\]
Then $S_n(e_1,\ldots,e_n)=n!e_1\cdots e_n$ and $\prod S_{r_i}(e_1,\ldots,e_{r_i})\not=0$ in $C$. The polynomial $\prod S_{r_i}(x_1,\ldots,x_{r_i})$
generates a submodule $N_m(\lambda)$ of $F_m$ and it does not vanish on $(C,V)$. Hence all submodules $N_m(\lambda)$ do participate in $F_m$ and $U$ coincides
with the $GL$-ideal of the weak identities for $(C,V)$. The proof of the theorem is completed.
\hfill $\Box$

{\it Proof of Corollary 1.} It follows from the proof of Theorem 1 that $F_m=\sum N_m(\lambda_1,\ldots,\lambda_m)$.
Each of the modules $N_m(\lambda_1,\ldots,\lambda_m)$ is generated by a polynomial
\[
f_{\lambda}(x_1,\ldots,x_r)=\prod S_{r_i}(x_1,\ldots,x_{r_i}),
\]
where $r=r_1,r_2,\ldots,r_p$ are the lengths of the columns of $[\lambda]$. Since $\dim V_k=k$, for $r>k$ we have that $S_r(x_1,\ldots,x_r)=0$
is a weak identity for the pair $(C_k,V_k)$. On the other hand $S_r(e_1,\ldots,e_r)\not=0$ in $C_r$ for $r\leq k$ and for every basis
$e_1,\ldots,e_k$ of $V_k$ with $\langle e_i,e_i\rangle\not=0$, $\langle e_i,e_j\rangle=0$, $i\not=j$. Therefore
\[
K\langle x_1,\ldots,x_m\rangle/(K\langle x_1,\ldots,x_m\rangle\cap T(C_k,V_k))=\sum N_m(\lambda_1,\ldots,\lambda_k)
\]
and the $GL$-ideal $T(C_k,V_k)$ is generated by $[x_1^2,x_2]$ and $S_{k+1}(x_1,\ldots,x_{k+1})$.\hfill $\Box$

{\it Proof of Corollary 2.} It is known \cite{R} that the pair $(M_2,sl_2)$ satisfies the weak identity (1). Besides (see \cite{DV})
\[
K\langle x_1,\ldots,x_m\rangle/(K\langle x_1,\ldots,x_m\rangle\cap T(M_2,sl_2)=\sum N_m(\lambda_1,\lambda_2,\lambda_3).
\]
The elements of the submodule $N_m(\lambda_1,\ldots,\lambda_m)$ of $F_m$ follow from $S_4(x_1,x_2,x_3,x_4)=0$ when $\lambda_4\not=0$.
This gives the assertion (i). The assertion (ii) follows immediately from (i) since one can easily obtain that the Lie weak identity
$S_4(x_1,x_2,x_3,x_4)=0$ is a consequence of (1). \hfill $\Box$


{\it Proof of Theorem 2.} It follows easily from (3) that the Young diagrams form a partially well-ordered set with respect to the inclusion.
Therefore each set $P$ of diagrams has a finite subset $P_0=\{[\lambda^{(1)}],\ldots,[\lambda^{(s)}]\}$ such that for every $[\lambda]\in P$
there exists a $[\lambda^{(i)}]\in P_0$ with $[\lambda^{(i)}]\subseteq [\lambda]$. Assume $I$ is a $GL$-ideal of $K\langle X\rangle$ containing $[x_1^2,x_2]$
and $\overline{I}$ is the homomorphic image of $I$ under the canonical homomorphism $K\langle X\rangle\to F=K\langle X\rangle/U$.
Let $N(\lambda)$ be the irreducible $GL$-module related to $\lambda$. It follows from the decomposition $F=\sum N(\lambda)$ that
$\overline{I}=\sum N(\lambda)$, where $[\lambda]$ ranges over some set $P$.

Denote by $P_0$ the finite subset of the elements of $P$ minimal with respect to the inclusion. The theorem will be proved if we establish that all
elements of $\overline{I}$ are consequences of $f_{\lambda}\in N(\lambda)$ where $[\lambda]\in P_0$.

As in \cite{K} it suffices to show that $f_{\mu}$ follows from $f_{\lambda}$ in $F$ when $[\mu]$ is obtained from $[\lambda]$ by adding a box. Let
$\lambda=(\lambda_1,\ldots,\lambda_r)$ and $m>r$. We define a homomorphism of $GL$-modules $\varphi:N_m(\lambda)\otimes N_m(1)\to F$
in the following way:
\[
\varphi\left(\sum f_i\otimes x_i\right)=\sum f_ix_i,\quad f_i\in N_m(\lambda).
\]
By the Branching Theorem $N_m(\lambda)\otimes N_m(1)=\sum N_m(\mu)$, where $[\mu]$ is obtained from $[\lambda]$ by adding a box.
The explicit form of the generators $f_{\mu}(x_1,\ldots,x_m)$ of $N_m(\mu)$ is found in \cite{K}. As in \cite{K} we can verify  that
$f_{\mu}(e_1,\ldots,e_m)\not=0$ for every basis $e_1,e_2,\ldots$ of $V$, $\langle e_i,e_i\rangle\not=0$, $\langle e_i,e_j\rangle=0$, $i\not=j$.
Thus we obtain that $f_{\mu}$ is a consequence of $f_{\lambda}$ in $F$. 
\hfill $\Box$

\end{document}